\documentclass {amsart}

\usepackage{setspace}

\usepackage{amsmath}
\usepackage{amssymb}
\usepackage{amsthm}
\usepackage{graphicx}

\input xy
\xyoption{all}

\newcommand{\g}[1]{\mathfrak{#1}}
\newcommand{\mc} [1]{\mathcal{#1}}

\newcommand{\pf} {\text{pf\ }}

\theoremstyle{plain}
\newtheorem{theorem}{Theorem}[section]
\newtheorem{lemma}[theorem]{Lemma}

\newtheorem{cor}[theorem]{Corollary}
\newtheorem{prop}[theorem]{Proposition}
\newtheorem{conjecture}[theorem]{Conjecture}

\theoremstyle{definition}
\newtheorem{definition}[theorem]{Definition}
\newtheorem{example}[theorem]{Example}

\newtheorem{remark}[theorem]{Remark}

\begin{document}

\title {On self-associated sets of points in small projective spaces}
\author {Ivan Petrakiev}
\address{Department of Mathematics, Harvard University, Cambridge, MA 02138}
\email{petrak@math.harvard.edu}
\begin {abstract} We study moduli of ``self-associated'' sets of points in ${\bf P}^n$ for small $n$. In particular, we show that for $n=5$ a general such set arises as a hyperplane section of the Lagrangean Grassmanian $LG(5,10) \subset {\bf P}^{15}$ (this was conjectured by Eisenbud-Popescu in {\it Geometry of the Gale transform}, J. Algebra 230); for $n=6$, a general such set arises as a hyperplane section of the Grassmanian $G(2,6) \subset {\bf P}^{14}$. We also make a conjecture for the next case $n=7$. Our results are analogues of Mukai's characterization of general canonically embedded curves in ${\bf P}^6$ and ${\bf P}^7$, resp.
\end{abstract}
\thanks{The author was partially supported by the NSF Graduate Research Fellowship.}

\maketitle

\section {Introduction}

Let $\Gamma \subset {\bf P}^n = {\bf P} V$ be a nondegenerate set of $2n+2$ distinct points with the following property: one can write $\Gamma = \Gamma' \cup \Gamma''$ such that both $\Gamma'$ and $\Gamma''$ correspond to orthogonal bases for some nondegenerate bilinear form $Q$ on $V$ (i.e. both $\Gamma'$ and $\Gamma''$ are {\it apolar} simplices with respect to $Q$). Then, we say that $\Gamma$ is {\it a self-associated (s.a.) set of points} in ${\bf P}^n$. This notion of self-association was first introduced by Castelnuovo \cite{Cast2} and further studied by Coble (\cite{Coble}, \cite{Coble2}) and other classical geometers (\cite{Bg},\cite{Bath}). The theory of self-association was given a modern treatment by Dolgachev-Ortland (\cite{Dolg}) and, more recently, by Eisenbud-Popescu  (\cite{EP}) in connection with the minimal resolution conjecture. Self-associated sets of points were also studied by Schreyer-Tonoli (\cite{ST}) in connection with the Green's conjecture on syzygies of canonical curves.

An equivalent definition of self-association is the following one (see \cite{EP}, thm. 7.1 and 8.1):

\begin{definition} Let $\Gamma \subset {\bf P}^{n}$ be a nondegenerate set of $2n+2$ distinct points in ${\bf P}^n$. Then, $\Gamma$ is s.a. set of points if and only if any subset of $2n+1$ points of $\Gamma$ imposes the same number of conditions on quadrics as $\Gamma$.
\end{definition}

There is a moduli space, which we denote ${\mc A}_n$, parametrizing (unordered) self-associated sets of points in ${\bf P}^n$ modulo projective equivalence. It is shown in \cite{Dolg}, that ${\mc A}_n$ is an irreducible, unirational variety of dimension $\binom {n+1} 2$.

A more sophisticated characterization of self-association is provided by the following theorem:

\begin{theorem}(\cite{EP}) Let $\Gamma$ be a nondegenerate set of $2n+2$ points in ${\bf P}^n$. Then, $\Gamma$ is Arithmetically Gorenstein (AG) if and only if $\Gamma$ is self-associated and $\Gamma$ fails exactly by 1 to impose independent conditions on quadrics.
\end{theorem}

It is easy to see, that if $\Gamma$ is a sufficiently general s.a. set of points (e.g. $\Gamma$ is in linearly general position), then the theorem above applies to $\Gamma$. 

\begin{example} If $\Gamma$ is a set of $2n+2$ points that lies on a rational normal curve $D \subset {\bf P}^n$, then $\Gamma$ is a s.a. set of points. 
\end{example}

\begin{example} Let $C \subset {\bf P}^{n+1}$ be a smooth canonically embedded curve. It is well-known, that $C$ is AG. Therefore, any transversal hyperplane section $\Gamma = C \cap {\bf P}^n$ is AG, and so, it forms a s.a. set of points in ${\bf P}^n$ (see \cite {ST}).
\end{example}

A geometric characterization of s.a. sets of points is classically known in small dimensions (see \cite{EP}, sec. 9). Namely, in ${\bf P}^2$, a s.a. set of points $\Gamma$ is just a complete intersection of a conic and a cubic. In ${\bf P}^3$, a sufficiently general s.a. set of points is a complete intersection of three quadrics. In ${\bf P}^4$, a sufficiently general s.a. set of points is a complete intersection of an elliptic normal curve $E \subset {\bf P}^4$ and a quadric. In particular, in all three cases, $\Gamma$ in arises as a hyperplane section of a canonical curve.
\bigskip

The goal of the present work is to give a similar characterization of s.a. sets of points in dimensions ${\bf P}^5$ and ${\bf P}^6$. More precisely, we will show, that a general s.a. set of points in these dimensions arises as a hyperplane section of a smooth canonical curve. In fact, our result for the case ${\bf P}^5$ answers a conjecture of Eisenbud-Popescu (\cite{EP}). On the other hand, the pattern ceases to generalize in dimension ${\bf P}^7$.
\bigskip

The starting point is the well-known theorem of Mukai (\cite{Mu2}):

\begin{theorem}\label{thm_muk} \begin{enumerate}
\item A general canonical curve $C \subset {\bf P}^6$ of genus $g=7$ is a linear section of the Lagrangean Grassmanian $LG_+(5,10) \subset {\bf P}^{15}$;
\item A genral canonical curve $C \subset {\bf P}^7$ of genus $g=8$ is a linear section of the Grassmanian $G (2,6) \subset {\bf P}^{14}$;
\item A general canonical curve $C \subset {\bf P}^8$ of genus $g=9$ is a linear section of the symplectic Grassmanian $G_\omega(3,6) \subset {\bf P}^{13}$.
\end{enumerate}
\end{theorem}

We will show that parts (a) and (b) of the theorem above extend to s.a. sets of points in ${\bf P}^5$ and ${\bf P}^6$, resp. On the other hand, the analogue of part (c) fails for s.a. sets of points in ${\bf P}^7$ because of a simple moduli count (see conjecture \ref{conj_sa}).
\bigskip

Ranestad-Schreyer (\cite{RS}) proved that ``a general empty AG scheme of degree 12 in ${\bf P}^4$'' (that is, a graded Artinian Gorenstein ring with Hilbert function (1,5,5,1)) arises as a linear section of $LG_+ (5,10)$. Thus, our result for s.a. sets of points in ${\bf P}^5$ should be viewed as an intermediate step between their result and that of Mukai. 

Also, we reinterpret our result for s.a. sets of points in ${\bf P}^5$ as an instance of ``inversion'' of the Cayley-Bacharach theorem for vector bundles (
%\cite {GrHa}, 
\cite{GH}, \cite{Laz}).

\subsubsection*{Notation and conventions.}

We work over a base field $k$ of characteristic $0$.
For any closed reduced subscheme $X \subset {\bf P}^n$ with ideal sheaf $\mc{I}_X$, and any positive integer $m$, denote by $X^{(m)} \subset {\bf P}^n$ the subscheme defined by the ideal sheaf $\mc{I}_X^{\otimes m}$. In other words, $X^{(m)}$ the smallest scheme supported on $X$ that contains the $(m-1)$-st infinitesimal neighborhood of $X$.
For any $r < n$, denote by $G (r,n) = G ({\bf P}^{r-1}, {\bf P}^{n-1})$ the Grassmanian of $(k-1)$-planes in ${\bf P}^{n-1}$.
\newpage

\section {Characterization of s.a. sets of points in small dimensions}

\subsection {The case ${\bf P}^5$}\

\bigskip

%\subsubsection* {The geometry of $LG_+ (5,10)$}

Recall some standard facts about Lagrangean Grassmanians (\cite{Mu}, \cite{FH}). Let $(V,Q)$ be a $2n$-dimensional vector space $V$ equipped with a nondegenerate symmetric bilinear form $Q$. A linear subspace $W \subset V$ is {\it isotropic} iff $Q(w,w) = 0$ for any $w \in W$. An isotropic subspace of maximal dimension $n$ is called a {\it Lagrangean} subspace. Since we assume the dimension of $V$ to be even, the variety parametrizing all such subspaces consists of two disjoint components $LG_\pm (n,V)$. 

Write $\wedge^{\bullet} V = \wedge^{ev} V \oplus \wedge^{odd} V$. Then, the Lagrangean Grassmanian $LG_+(n,V)$ (resp. $LG_-(n,V)$) admits the Pl\"ucker embedding $LG_+ (n,V) \hookrightarrow {\bf P} (\wedge^{ev} V)$ (resp. $LG_- (n,V) \hookrightarrow {\bf P} (\wedge^{odd} V)$). The special orthogonal group $SO(V,Q)$ acts on $\wedge^{ev} V$ (resp. $\wedge^{odd} V$) via the {\it (half-)spin representation} and $LG_+ (n,V)$ (resp. $LG_- (n,V)$) is a homogeneous space for this representation. 

Next, we specialize to the case $\dim V = 10$. Denote $X = LG_+ (5,10)$. One can  show (\cite{Mu}), that $X \subset {\bf P}^{15}$ is a 10-dimensional variety of degree 12. Also, the canonical bundle $K_X \cong \mc{O}_X (-8)$. 

It is observed in (\cite{Mu}), that the homogeneous ideal of $X \subset {\bf P}^{15}$ is generated by $10$ quadrics, that in turn satisfy a unique quadratic relation. More precisely, consider the natural map
$$\rho: Sym^2 H^0 ({\bf P}^{15}, \mc{I}_X (2)) \longrightarrow H^0 ({\bf P}^{15}, \mc{I}_{X^{(2)}}(4))$$
induced by multiplication. Then, the kernel of $\rho$ is 1-dimensional, generated by a nondegenerate symmetric bilinear form, which we denote $R^*$. Taking the dual, the space of quadrics $H^0 ({\bf P}^{15}, \mc{I}_X(2))$ is naturally endowed with a nondegenerate symmetric bilinear form $R$.

Notice, that there is a natural action of $SO(10)$ on $H^0 ({\bf P}^{15}, \mc{I}_X(2))$, induced by the spin representation.

\begin{prop}\label{prop_mukai}(Mukai, \cite{Mu}) There is a natural isomorphism
$$\alpha : V \stackrel {\sim}{\longrightarrow} H^0 ({\bf P}^{15}, \mc{I}_X(2))$$
with the following properties:
\begin{enumerate}
\item $\alpha$ is $SO(10)$-equivariant and $\alpha^* (R) = Q$.
\item Let $x \in X$ be a point, corresponding to a Lagrangean $W_x \subset V$. Then, $\alpha(W_x) \subset H^0 ({\bf P}^{15}, \mc{I}_X(2))$ is precisely the linear subsystem of quadrics containing $X$ that are singular at $x$.
\end{enumerate}
\end{prop}

%\subsubsection* {Self-associated sets of points in ${\bf P}^5$} 

Let $X = LG_+ (5,10)$ as above. Since $K_X \cong \mc{O}_X(-8)$ and $X$ is projectively normal, it follows, that a general linear section $\Gamma = X \cap {\bf P}^5$ is a s.a. set of points in ${\bf P}^5$. It was conjectured by Eisenbud-Popescu in \cite {EP}, that any sufficiently general s.a. set of points in ${\bf P}^5$ arises in this way. Indeed, we have:

\begin{theorem}\label{thm_12} Let $\Gamma$ be a general s.a. set of points in ${\bf P}^{5}$. Then, $\Gamma$ is a linear section of the Lagrangean Grassmanian $LG_+ (5,10) \subset {\bf P}^{15}$. Moreover, the linear section is unique upto the natural $SO(10)$-action induced by the spin representation.
\end{theorem}

From Mukai's result, we conclude:
\begin{cor}
A general s.a. set of points in ${\bf P}^5$ is a linear section of a canonical curve in ${\bf P}^6$.
\end{cor}

Our approach is essentially a translation of Mukai's proof from curves to points. A certain geometric argument of Mukai is replaced by a symbolic computation on the computer algebra system {\it Macaulay 2} (\cite {Mac2}). 
\bigskip

{\it Proof of Theorem \ref{thm_12}.} Let $\Xi \subset G({\bf P}^5, {\bf P}^{15})$ be the open subset of linear subspaces that meet $X \subset {\bf P}^{15}$ transversally at 12 distinct points. Notice, that $SO(10)$ acts on $\Xi$ and there is a well-defined morphism
$$\sigma : \Xi / SO (10) \longrightarrow \mc{A}_5,$$
that assigns to every $[\Lambda] \in \Xi$ the s.a. set of points $X \cap \Lambda$.
By counting moduli, we get: $$\dim \Xi / SO(10) = 6(15-5) - \binom {10} 2 = 15$$
and
$$\dim \mc{A}_5 = \binom{5+1} 2 = 15.$$
Therefore, to show that $\sigma$ is birational, it suffices to show that $\sigma$ is injective. 

We will need the following result.

\begin{lemma}\label{lemma_horrible} The Lagrangean Grassmanian $X \subset {\bf P}^{15}$ is ACM and so is the scheme $X^{(2)} \subset {\bf P}^{15}$.
\end{lemma}
\begin{proof} Homogeneous varieties are known to be ACM in general. We have verified that $X^{(2)}$ is ACM by using {\it Macaulay 2}. See the Appendix for details.
\end{proof}

\begin{cor}\label{cor_commute} Let $\Gamma = X \cap {\bf P}^5$ be any transversal linear section. Consider the commutative diagram
$$\xymatrix{
Sym^2 H^0 ({\bf P}^{15}, \mc{I}_X (2)) \ar[r]^\rho \ar[d] & H^0 ({\bf P}^{15}, \mc{I}_{X^{(2)}}(4)) \ar[d] \\
 Sym^2 H^0 ({\bf P}^{5}, \mc{I}_\Gamma(2)) \ar[r]^{\overline \rho} &  H^0 ({\bf P}^{5}, \mc{I}_{\Gamma^{(2)}} (4))\\
}
$$
where the horizontal maps are induced by multiplication and the vertical maps are induced by restriction. Then, the two vertical maps are isomorphisms.
\hfill$\Box$
\end{cor}

We are ready to complete the proof of the theorem. Let $\Gamma = X\cap \Lambda \subset {\bf P}^{15}$ be a transversal linear section with $\Lambda \cong {\bf P}^5$. By cor. \ref{cor_commute}, the map
$$\overline {\rho}: Sym^2 H^0 ({\bf P}^{5}, \mc{I}_\Gamma(2)) \rightarrow H^0 ({\bf P}^{5}, \mc{I}_{\Gamma^{(2)}} (4))$$
has a 1-dimensional kernel, generated by a nondegenerate symmetric bilinear form ${\bar R}^*$. Let ${\bar R}$ be the dual of ${\bar R}^*$. 

[Here is a heuristic way to see why $\dim \ker \overline \rho = 1$. Clearly, $\dim Sym^2 H^0 ({\bf P}^5, \mc{I}_\Gamma(2)) = 55$. Also, we would have $h^0 ({\bf P}^5, \mc{I}_{\Gamma^{(2)}} (4)) = \binom {4+5}{5} - 12 \cdot 6 = 54$, assuming that double points at $\Gamma$ impose independent conditions on quartics in ${\bf P}^5$.]

Let
$$\beta : H^0 ({\bf P}^{15}, \mc{I}_X(2)) \stackrel {\sim}{\longrightarrow} H^0 ({\bf P}^{5}, \mc{I}_\Gamma(2))$$
be the map, induced by restriction to $\Lambda$. By cor. \ref{cor_commute}, we have $\beta^*(\bar R) = R$. Composing with the map $\alpha$ from prop. \ref{prop_mukai}, we get an isomorphism
$$\beta \circ \alpha: V \stackrel {\sim}{\longrightarrow} H^0 ({\bf P}^{5}, \mc{I}_\Gamma(2))$$
such that $(\beta\circ \alpha)^* (\overline R) = Q$. 

From now on, we identify $V$ with $H^0 ({\bf P}^5, \mc{I}_\Gamma(2))$ via $\beta \circ \alpha$. For any $p_i \in \Gamma$, denote by $W_i \subset V$ the linear subspace parametrizing quadrics containing $\Gamma$ that are singular at $p_i$. Since $\Gamma$ is cut-out scheme-theoretically by quadrics, a double point at $p_i \in \Gamma$ imposes $5$ conditions on quadrics in $\Lambda$ passing through $p_i$. Therefore, $\dim W_i = 5$. It follows, that any such quadric in $\Lambda$ is a restriction of a quadric in ${\bf P}^{15}$ containing $X$ and singular at $p_i$. By prop. \ref{prop_mukai}, $W_i \subset V$ is a Lagrangean and moreover, it corresponds to the point $p_i \in \Gamma \subset X$, i.e. $p_i = [W_i]$.

%\begin{center}
%\ \ \ \ \ \ \ \ \includegraphics[0cm,0cm][6cm,6cm]{spinor.eps}
%\end{center}
%\begin{center}
%{\bf Figure 3}
%\end{center}

Now, let $\Gamma' = X \cap \Lambda'$ be another transversal linear section and assume, that $\Gamma \subset \Lambda$ and $\Gamma' \subset \Lambda'$ are projectively equivalent. We claim, that $\Lambda$ and $\Lambda'$ are in the same $SO(10)$-orbit.

For, let
$$\beta' : H^0 ({\bf P}^{15}, \mc{I}_X(2)) \stackrel {\sim}{\longrightarrow} H^0 ({\bf P}^{5}, \mc{I}_{\Gamma'}(2))$$
be the map, induced by restriction to $\Lambda'$. Hence, $\beta^* (\bar R) = \beta'^* (\bar R) = R$. It follows, that $\beta' = \beta \circ \tau$ for some $\tau \in SO(10)$. For any $p_i' \in \Gamma' \subset X$, define the Lagrangean $W'_i \subset V$ as before. It follows, that $[W'_i] = \tau [W_i]$. Therefore, $[\Lambda'] = \tau [\Lambda]$.

This shows, that the map $\sigma: \Xi / SO(10) \rightarrow {\mc A}_5$ is injective, which completes the proof. \hfill$\Box$

\bigskip

\subsubsection*{Alternative interpretation}\

\bigskip

We interpret Theorem \ref{thm_12} as an instance of ``inversion'' of the Cayley-Bacharach theorem for vector bundles (see 
%\cite{GrHa}, 
\cite{GH}, \cite{Laz}). 

Consider a surjective map vector bundles
$$f : \mc{O}_{{\bf P}^5} (2)^7 \rightarrow \mc{O}_{{\bf P}^5}(3)^2$$
and let $B_f = \ker f$. If $f$ is  sufficiently general, then $B_f$ is a vector bundle of rank 5 on ${\bf P}^5$. By Whitney's product formula, we have:
\begin{align*}
c_1 (B_f) &= 8[H]; \\ 
c_5 (B_f) &= 12 [H]^5.
\end{align*}Here, $[H]$ denotes the hyperplane class.

Let $\Gamma = Z(s)$ be the zero-locus of a regular section $s \in H^0 (B_f)$ and assume that $\Gamma$ is reduced. By the Cayley-Bacharach theorem for vector bundles, $\Gamma$ fails exactly by one to impose independent conditions on the linear system $|K_{{\bf P}^5} \otimes \det(B_f)| = |\mc{O}_{{\bf P}^5} (2)|$. In other words, $\Gamma$ is a self-associated set of points in ${\bf P}^5$.

Remarkably, Theorem \ref{thm_12} implies (and, in fact, is equivalent to) the converse statement: 
\smallskip

{\it If $\Gamma$ is a general self-associated set of points in ${\bf P}^5$, then $\Gamma$ is the zero locus of a regular section of a vector bundle $B_f$, for a suitable choice of $f$.}
\bigskip
 
We just sketch the proof. Consider a sufficiently general map of vector bundles 
$$\widetilde f : \mc{O}_{{\bf P}^{15}} (2)^7 \rightarrow \mc{O}_{{\bf P}^{15}}(3)^2$$
and let $B_{\widetilde f} = \ker \widetilde f$. Then, $B_{\widetilde f}$ is a reflexive sheaf of rank 5 (in fact, a Buchsbaum-Rim sheaf, in the language of \cite{Mig}). Let $s \in H^0 (B_{\widetilde f})$ be a general section. One can show, that the top component of the zero locus $Z(s)$ is projectively equivalent to the Lagrangean Grassmanian $LG_+ (5,10)$ ($Z(s)$ also has a component of larger codimension, namely the locus of degeneracy of $\widetilde f$). We have checked this on {\it Macaulay 2}. It remains to consider the restriction of $B_{\widetilde f}$ to a general ${\bf P}^5 \subset {\bf P}^{15}$ and apply Theorem \ref{thm_12}.

\subsection {The case ${\bf P}^6$} \

\bigskip

Recall some standard facts about Grassmanians of lines (\cite{HFirst}). Let $V$ be an $n$-dimensional vector space over the base field $k$. Consider the Grassmanian $X = G(2,V) \subset {\bf P} (\wedge^2 V)$ parametrizing lines in ${\bf P} V$, together with its Pl\"ucker embedding. The group $SL(V)$ naturally acts on $\wedge^2 V$ and $X$ is a homogeneous space for this action.

For any $w\in \wedge^2 V$ and $k \leq n/2$, denote $$\pf (2k, w) = \mathop {\underbrace {w \wedge \dots \wedge w}}_{k \ \text{times}} \in \wedge^{2k} V.$$
If we fix a basis for $V$, then $\pf(2k,w)$ is just the $2k\times 2k$-pfaffians of the skew-symmetric matrix given by $w$. 

It is well-known, that $X \subset {\bf P} (\wedge^2 V)$ can be defined as the locus of vanishing of $\pf(4,w)$, for $w \in \wedge^2 V$. More generally, for any $m$, the $m$-secant variety $\sigma_m (X)$ can be defined as the locus of vanishing of $\pf(2m+4, w)$. 

Next, we specialize to the case $\dim V = 6$. Then, $X = G(2,6) \subset {\bf P}^{14}$ is an 8-dimensional variety of degree 14. Also, the canonical bundle $K_X \cong \mc{O}_X(-6)$. It follows, that a general section $\Gamma = X \cap {\bf P}^6$ is a s.a. set of points in ${\bf P}^6$. Our main result in this section is the converse statement:

\begin{theorem}\label{thm_P6} A general s.a. set of points in ${\bf P}^{6}$ is a linear section of the Grassmanian $G(2,6) \subset {\bf P}^{14}$.
\end{theorem}

By Mukai's theorem, we have:

\begin{cor} A general s.a. set of points in ${\bf P}^6$ is a hyperplane section of a canonical curve in ${\bf P}^7$. 
\end{cor}

{\it Proof of theorem \ref{thm_P6}}. Just as in the proof of theorem \ref{thm_12}, we start with moduli count. Let $\Xi \subset G ({\bf P}^6, {\bf P}^{14})$ be the open subset of linear subspaces that meet $X \subset {\bf P}^{14}$ transversally at $14$ distinct points. Now, $SL(6)$ acts on $\Xi$ and we have a well-defined map
$$\sigma : \Xi / SL(6) \rightarrow \mc{A}_6,$$
that assigns to every $[\Lambda] \in \Xi$ the s.a. set of points $X \cap \Lambda$. We have: $$\dim \Xi / SL(6) = 7(14-6) - 35 = 21$$ and $$\dim \mc{A}_6 = \binom{6+1} 2 = 21.$$
To complete the proof of the theorem, we will show, that $\sigma$ is generically finite. 

We will need the following result.

\begin{lemma} The Grassmanian $X = G(2,6) \subset {\bf P}^{14}$ is ACM and so is the scheme $X^{(2)}$.
\end{lemma}

\begin{proof} Since $X$ is homogeneous, it is ACM. We have verified the statement for $X^{(2)}$ via {\it Macaulay 2} (See the Appendix for details). 
\end{proof}

Consider the secant variety $\sigma_1 (X)$, given by $pf(6,w) = 0$. Then, $\sigma_1(X)$ is the unique cubic hypersurface in ${\bf P}^{14}$ with singular locus $X$. Since $X^{(2)}$ is ACM, we get:

\begin{cor}\label{pff} Let $\Gamma = X \cap {\bf P}^6$ is any transversal linear section. Then, there is a unique cubic hypersurface in ${\bf P}^6$ that is singular at every point of $\Gamma$.
\hfill$\Box$
\end{cor}

Let $\Gamma = X \cap \Lambda$ for a sufficiently general $[\Lambda] \in \Xi$. Then, $\Gamma$ is defined as the vanishing locus of the $4\times 4$ pfaffians in a $6\times 6$ skew-symmetric matrix $A \in H^0 (\wedge^2 V \otimes \mc{O}_{{\bf P}^6}(1))$. In coordinate-free form, $A$ is simply the restriction of the identity element $1 \in H^0 (\wedge^2 V \otimes \mc{O}_{{\bf P}(\wedge^2 V)} (1))$. Conversely, the matrix $A$ uniquely determines $[\Lambda] \in \Xi$.

Now, let $\Gamma' = X \cap \Lambda'$ for some other $[\Lambda'] \in \Xi$ and assume, that $\Gamma' \subset \Lambda'$ is projectively equivalent to $\Gamma \subset \Lambda$. Let $B \in H^0 (\wedge^2 V \otimes \mc{O}_{{\bf P}^6}(1))$ be the corresponding skew-symmetric matrix. In what follows, we identify the ambient spaces of $\Gamma$ and $\Gamma'$, and simply write ${\bf P}^6$. After normalizing $B$ by a scalar if necessary, cor. \ref{pff} implies, that $pf(6,A) = pf(6,B)$. We have:

\begin{lemma}\label{lemma_local} Let $\Gamma$ and $A \in H^0 (\wedge^2 V \otimes {\mc O}_{{\bf P}^6}(1))$ be as above. Then, there is a Zariski open subset $U \subseteq H^0 (\wedge^2 V \otimes {\mc O}_{{\bf P}^6}(1))$, such that the following is true: for any $B \in U$, $\pf(6,A) = \pf(6,B)$ if and only if $A = B^\tau$ for some $\tau \in SL(V)$.
\end{lemma}

(Here we denote $A^\tau := \tau^t A \tau$, which is independent of the choice of a basis for $V$.)

\begin{proof} Since $\pf(6, A) = \pf(6, A^\tau)$, the ``if'' part is true for any open $U$.

Let $Z \subset H^0 (\wedge^2 V \otimes {\mc O}_{{\bf P}^6}(1))$ be the closed subset of elements $B$ such that $\pf(6,A) = \pf(6,B)$, and let $T_A Z$ be the Zariski tangent space to $Z$ at $A$. Consider the map $\rho: SL(V) \rightarrow Z$, given by $\tau \mapsto A^\tau$, and its differential $d\rho: \g{sl}(V) \rightarrow T_A Z$. It is easy to see, that if $A$ is sufficiently general, then $d\rho$ is injective.

Next, we will show, that $d\rho$ is in fact bijective. For any $B' \in T_A Z$, let $B = A + \epsilon B'$, where $\epsilon^2 = 0$. Since $\pf(B) = \pf(A) + 3 \epsilon B' \wedge A \wedge A$ and $\text{char}\ k \neq 3$, we have:
$$B' \wedge A \wedge A = 0.$$
This is a system of linear equations for $B'$ which can be solved explicitly by looking at the resolution of $\Gamma$. Since $X$ is ACM, $\Gamma$ has the same resolution as $X$ (see the Appendix):
$$ \dots \longrightarrow \mc{O}_{{\bf P}^6}(-3)^{35} \longrightarrow \mc{O}_{{\bf P}^6}(-2)^{15} \stackrel {A \wedge A}{\longrightarrow} \mc{O}_{{\bf P}^6} \longrightarrow \mc{O}_\Gamma \longrightarrow 0$$
Therefore, $\dim T_A Z = \dim \g{sl} (V) = 35$, and so $d \rho : \g{sl}(V) \rightarrow T_A Z$ is bijective.

The lemma follows.
\end{proof}

From what we said, it follows that there is an open subset $\widetilde U \subset \Xi$ containing $[\Lambda]$ such that the following is true: if $[\Lambda'] \in \widetilde {U}$ is such that $\Gamma' = X \cap \Lambda'$ is projectively equivalent to $\Gamma$, then $[\Lambda'] = \tau [\Lambda]$, for some $\tau \in SL(V)$. 

This implies, that the map $\sigma : \Xi / SL(V) \rightarrow \mc{A}_6$ is generically finite, which completes the proof of the theorem. \hfill$\Box$

\begin{remark} We expect, that the map $\sigma : \Xi /SL(V) \rightarrow \mc{A}_6$ is actually birational. This would be true, if we knew that lemma \ref{lemma_local} holds globally, i.e. we may take $U = \Gamma(\wedge^2 V \otimes {\mc O}_{{\bf P}^6}(1))$. This seems likely, but we don't know how to prove it.
\end{remark}

%Finally, in the spirit of Ranestad-Schreyer (\cite {RS}), we ask:
%
%\begin{question} Does ``a general empty AG scheme of degree 14 in ${\bf P}^5$'' (that is, a graded Artinian Gorenstein ring with Hilbert function (1,6,6,1)) arise as a linear section of $G(2,6) \subset {\bf P}^{14}$?
%\end{question}

\subsection {The case ${\bf P}^7$}\

\bigskip

As we saw, for any $n \leq 6$, a general s.a. set of points in ${\bf P}^n$ is a hyperplane section of a canonical curve of genus $n+2$ in ${\bf P}^{n+1}$. Clearly, this pattern cannot continue for $n$ arbitrary large, simply because $\dim \mc{M}_g = 3g-3$ and $\dim \mc{A}_{g-2} = \binom{g-1} 2$. In fact, $n=7$ is the first case when a general s.a. set of points in ${\bf P}^7$ is {\it not} a section of a canonical curve.

To see this, recall Mukai's result that a general canonical curve of genus $9$ in ${\bf P}^8$ is a linear section of the isotropic Grassmanian $G_\omega (3,6) \subset {\bf P}^{13}$. We claim, that a general s.a. set of points $\Gamma \subset {\bf P}^7$ is {\it not} a section of $G_\omega(3,6)$. For, let $\Xi \subset G({\bf P}^7, {\bf P}^{13})$ be the open subset of linear sections that meet $G_\omega (3,6)$ transversally at 16 distinct points. As before, consider the map $$\sigma : \Xi / Sp(6) \longrightarrow {\mc A}_7.$$ By counting moduli, we get $$\dim \Xi / Sp(6) = 8(13-7) - \binom{6+1} 2 = 27,$$ while $$\dim {\mc A}_7 = \binom {7+1} 2 = 28.$$ Therefore, $\sigma$ cannot be dominant!

Naturally, we make
\begin{conjecture}\label{conj_sa} $\sigma$ is generically injective.
\end{conjecture} 
If the conjecture is true, it would imply that there is a divisor $\mc{Y} \subset {\mc A}_7$ parametrizing those s.a. sets of points in ${\bf P}^7$ that arise as a hyperplane sections of canonical curves in ${\bf P}^8$. 

\begin{remark} We note a parallel with the question of whether a given canonical curve $C \subset {\bf P}^{g-1}$ arises as a hyperplane section of a polarized K3 surface (\cite{CiHaMi},\cite{Mu2},\cite{W}). It is well-known, that if $g \leq 11$, $g \neq 10$, then any canonical curve has this property. On the other hand, if $g=10$, then the curves with this property are parametrized by a divisor $\mc{K} \subset {\mc{M}_{10}}$ (\cite {FP}).
\end{remark}

\newpage
\subsection*{APPENDIX: Some computations on Macaulay 2}\

\bigskip

Below are listed the Betti numbers of $LG_+ (5,10) \subset {\bf P}^{15}$ and $G(2,6) \subset {\bf P}^{14}$ and the schemes of their first infinitesimal neighborhoods. In particular, we see that all schemes in consideration are ACM.

\begin{singlespace}
$$
%\footnotesize
\begin{array}{c|cccccc}
\text{degree} & \\
\hline
0 & 1 & - & - & - & - & - \\
1 & - & 10 &16 & - & - & - \\
2 & - & - & - & 16 & 10 & - \\
3 & - & - & - & - & - & 1 \\
\end{array}
\ \ \ \ \ \ \ \ 
\begin{array}{c|ccccccc}
\text{degree} & \\
\hline
0 & 1 & - & - & - & - & - & - \\
1 & - & 15 & 35 & 21 & - & - & - \\
2 & - & - & - & 21 & 35 & 15 & - \\
3 & - & - & - & - & - & - & 1 \\
\end{array}
$$ 

$$LG_+ (5,10) \ \ \ \ \ \ \ \ \ \ \ \ \ \ \ \ \ \ \ \ \ \ \ \ \ \ \ \ \ \ \ \ \ \ \ \ \ G(2,6)$$
\bigskip

$$
%\footnotesize
\begin{array}{c|cccccc}
\text{degree} & \\
\hline
0 & 1 & - & - & - & - & - \\
1 & - & - & - & - & - & - \\
2 & - & - & - & - & - & - \\
3 & - & 54 & 144 & 120 & - & - \\
4 & - & - & - & - & 45 & 16 \\
\end{array}
\ \ \ \ \ \ \ \
\begin{array}{c|ccccccc}
\text{degree} & \\
\hline
0 & 1 & - & - & - & - & - & - \\
1 & - & - & - & - & - & - & - \\
2 & - & 1 & - & - & - & - & - \\
3 & - & 105 & 399 & 595 & 405 & 105 & - \\
4 & - & - & - & - & 21 & 35 & 15 \\
\end{array}
$$

$$LG_+ (5,10)^{(2)} \ \ \ \ \ \ \ \ \ \ \ \ \ \ \ \ \ \ \ \ \ \ \ \ \ \ \ \ \ \ \ \ \ \ \ \ \ G(2,6)^{(2)}$$
\end{singlespace}

\newpage

%--------------------------------------------------
%--------------------------------------------------
%--------------------------------------------------
%--------------------------------------------------

\bibliographystyle{amsplain}

\end{document}